\input amstex
\documentstyle{amsppt}
\magnification=\magstep1

\hsize=5.9 true in
\vsize8.8true in
\NoBlackBoxes
\define\GL{\operatorname{GL}}
\define\Alg{\operatorname{alg}}
\define\alg{\operatorname{alg}}
\redefine\span{\operatorname{span}}
\define\tr{\operatorname{tr}}
\define\PNC{\operatorname{PNC}}
\define\Log{\operatorname{Log}}
\define\sd{\operatorname{sd}}
\widestnumber\key{999}
\rightheadtext { }
\leftheadtext { }
\topmatter 
\title 
 On invertibility preserving linear mappings,\\
simultaneous triangularization and Property $L$
\endtitle 
\author 
Erik Christensen 
\endauthor 
\affil
Department of Mathematics
University of Copenhagen
Universitetsparken 5
DK 2100 Copenhagen
\endaffil
\endtopmatter 

\document

\subhead{1. Introduction}\endsubhead

The investigation leading to this publication was motivated by a desire to try
to understand the structure of a linear unital mapping $\varphi$ from a unital
algebra $\Cal A$ of matrices contained in $M_h(\Bbb C)$ into $M_n(\Bbb C)$
which has the property that an invertible element in $\Cal A$ is mapped into
an invertible in $M_n(\Bbb C)$. The interest in this question was raised by
some earlier results on a linear invertibility preserving mapping from a
Banach algebra into $M_n(\Bbb C)$. This lead to the article [3] which on the
other hand basically is of finite dimensional nature. During the search for
results of this type in the literature on linear algebra we found similar
questions discussed either as operators on $M_n(\Bbb C)$ [13], or in the form
of results on the structure of sets of matrices having the Property $L$ of
Motzkin and Taussky [10,11,12,15,16,17,18,19]. In [13] Marcus and Purves
describe the structure of an operator $\varphi$ on $M_n(\Bbb C)$ which
preserves invertibility. It turns out that $\varphi$ must have one of the
forms $$\varphi(A)=UAV\;\text{ or }\; \varphi(A)=UA^tV $$ where $U$ and $V$
are in $\GL_n(\Bbb C)$. If further $\varphi$ is unital $(\varphi(I)=I)$ then
$V=U^{-1}$ so one gets in this case that $\varphi$ is either an automorphism
or an anti-automorphism of the algebra $M_n(\Bbb C)$. 

Kaplansky [9] discussed the problem whether an invertibility preserving,
unital, continuous linear mapping between Banach algebras might have some
algebraic properties, and he suggested mainly on the basis of results from [5,
6, 13] that such mappings might be Jordan homomorphisms
$(\varphi(a^2)=\varphi(a)^2)$. Some experiments with low dimensional matrices
(Example 3.1) show that this is too much to expect. At first we thought that
Kaplansky was right modulo the Jacobson radical in the algebra generated by
the image of such a mapping. On the other hand a closer look at the simple
case where $\varphi$ is a mapping of $\Bbb C^3$ into $M_n(\Bbb C)$ shows that
this is not so. An example by Wielandt [19], which also has been used by
Motzkin and Taussky to give an $L$-pair which is not simultaneously
triangularizable yields an example which kills this modified conjecture
definitely. But still we had some results from [3] and some evidence from
notably [10,11,16,17] that some reasonable extra conditions might imply that
an invertibility preserving mapping with these extra properties is a Jordan
homomorphism.

Regarding this question we consider in section 4 a unital algebra $\Cal
A\subseteq M_h(\Bbb C)$ and a linear mapping $\varphi:\Cal A\to M_n(\Bbb C)$
such that $\varphi(I_{\Cal A})=I_{M_n(\Bbb C)}$ and the image of any
invertible element in $\Cal A$ is invertible in $M_n(\Bbb C)$. For $k$ in
$\Bbb N$ we define $\varphi_k:M_k(\Cal A)=\Cal A\otimes M_k(\Bbb C)\to
M_k(M_n(\Bbb C))=M_n(\Bbb C)\otimes M_k(\Bbb C)$ by $\varphi_k=\varphi\otimes
\text{id}$, i.e\.  $\varphi$ operates on each of the entries in $M_k(\Cal
A)$. We will then say that $\varphi$ is $k$-invertibility preserving if
$\varphi_k$ is invertibility preserving and our main result says that
$\varphi$ is a homomorphism modulo the Jacobson radical if and only if it is
$k$-invertibility preserving for some sufficiently large $k$. In particular
$\varphi$ will be $m$-invertibility preserving for any $m$ larger than this
$k$. The main problem is then to find an estimate for the least k which has
the property that $k$-invertibility implies that $\varphi$ is a homomorphism
modulo the radical. The bound becomes especially nice when the image is an
algebra. In other cases the bound, we have, depends on how far the image is
from being an algebra. We prove that $\varphi$ is a homomorphism if it is
$k$-invertibility preserving for $k=\dim(\alg(\varphi(\Cal
A))-\dim\varphi(\Cal A)+3$.  This estimate is fairly crude, so we introduce a
concept, which we call the semi-simple defect of a set of matrices, as a
measure for how far this set is from being an algebra modulo the Jacobson
radical. If $\varphi(\Cal A)=\alg(\varphi(\Cal A))$ then the semi-simple defect
is $0$ and $\varphi$ is a Jordan-homomorphism modulo the Jacobson radical
without further conditions, if it further is $2$-invertibility preserving then
it is a homomorphism modulo the Jacobson radical.
  
The Property $L$ was mentioned above. It is closely related to the
invertibility preserving question for an abelian algebra $\Cal A$.  In section
3 we introduce for $k\in\Bbb N$ a Property $kL$ which is stronger than $L$ and
mimics the ideas presented above in the way that Property $kL$ is a condition
on the characteristic roots for matrices of the form $a\otimes x+b\otimes y$
with $a,b$ fixed in $M_n(\Bbb C)$ and $x,y$ general in $M_k(\Bbb C)$. One can
say it is a set of conditions on matrix pencils [5, Matrizenb\"uschel] but
with coefficients in $M_k(\Bbb C)$ rather than in $\Bbb C$. We prove that for
some set $\Cal S$ contained in $M_n(\Bbb C)$ this set can be triangularized
simultaneously if and only if it has property $kL$ where $k=\dim(\alg(\Cal
S))-\dim(\span(\Cal S))+3$. As for the $k$-invertibility case it is relatively
easy to see that $(k+1)L$ implies $kL$ so the problem is to determine some
estimate for the lower bound for the set of those $k$ for which property $kL$
implies simultaneous triangularization. The least usable $k$ we find is the
semi-simple defect plus 3, but we find it is likely that it could be of the
order of the square root of the semi-simple defect rather than of first order.
 
  Finally section 2 contains the results which explain why we have such
estimates on $k$ and why it works. The content of Section 2 is not new, but
this section has its focus on traces, on the Jacobson radical and on algebras
of matrices rather than on a single matrix or pairs of matrices, and we have
not found this point of view presented in the literature on matrix theory.
Generally speaking one can say that section 2 repeats some well known
elementary algebra and tells you how much you have to pay in order to use
these results if you just have a set of matrices instead of an algebra. We
estimate the size of the degree of polynomials involved in order to get from a
subspace to the algebra generated by this subspace. In the sections 3 and 4
this estimate is used as the size $k$ on the matrix algebra we have to tensor
with in order to get the desired properties. The estimate given in section 2
is probably not optimal, on the other hand some rather concrete computations
made at the end of section 2 show that it may be hard to be more precise.

\subhead{2. Basic observations}\endsubhead

Through this section $\Cal B$ will denote an algebra of $n\times n$
matrices over $\Bbb C$ such that the unit $I$ of $M_n(\Bbb C)$ is in
$\Cal B$. By Wedderburns Theorem [2, p\. 143] and Wedderburn-Artins
Theorem [2, p\. 69], $\Cal B$ decomposes as a direct sum of a
semi-simple sub-algebra $\Cal C$ and the Jacobson radical $J$. Further
$\Cal C$ is isomorphic to $\Cal B/J$ and is a direct sum of full
matrix algebras hence $\Cal C=M_{n_1}(\Bbb C)\oplus\cdots\oplus
M_{n_s}(\Bbb C)$ and
 $$\Cal B=M_{n_1}(\Bbb C)\oplus\cdots\oplus M_{n_s}(\Bbb C)\oplus\Cal
J\,,\;\forall j\in\Cal J:\quad j^n=0\,.\tag1
 $$
 
It is important and will be used in the following arguments that in the sum
decomposition chosen $\Cal B=\Cal C\oplus\Cal J$, the algebra $\Cal C$ is a
sub-algebra of $\Cal B$. We will let $\tr_n$ denote the trace on $M_n(\Bbb C)$
and $I$ the unit in $M_n(\Bbb C)$ so $\tr_n(I)=n$. The restriction of $\tr_n$
to $\Cal J$ vanishes since the elements in $\Cal J$ are all nilpotent. The
restriction of $\tr_n$ to the summand $M_{n_i}(\Bbb C)$ in $\Cal C$ is a trace
on this algebra $\left(\tr_n(xy)=\tr_n(yx)\right)$ and for a minimal
idempotent $e$ in $M_{n_i}(\Bbb C)$ we have $d_i=\tr_n(e)\in\Bbb N$.  Since
all traces on $M_{n_i}(\Bbb C)$ are proportional to the canonical one --
$\tr_{n_i}$ -- we get $\tr_n|M_{n_i}(\Bbb C)=d_i\tr_{n_i}$. With this notation
in mind we can formulate the first observation, which is a known fact to
which we do not have a suitable reference.

\proclaim{2.1 Lemma} An element $b$ in $\Cal B$ belongs to $\Cal J$
if and only if for each $x\in\Cal B$ $\tr_n(xb)=0$.
\endproclaim

\demo{Proof} If $b$ is in $\Cal J$ then so is $xb$ and since $\Cal J$ is an
ideal,  $\tr_n$ vanishes on $xb$.

If $b$ does not belong to $\Cal J$ then according to the decomposition
of $\Cal B$ in (1)
 $$b=b_1+\cdots+b_i+\cdots+b_s+j
 $$ and we may assume that $b_i\neq0$. Let $$x=0+\cdots+0+b_i^\ast
+0+\cdots+0+0\,,
 $$
 where $b_i^\ast $ is the adjoint matrix in $M_{n_i}(\Bbb C)$ to $b_i$.
Then $x$ is in $\Cal B$ and
 $$\tr_n(xb)=d_i\tr_{n_i}(b_i^\ast b_i)>0\,.
 $$
\hfill$\square$
 \enddemo
 
\proclaim{2.2 Definition} For natural numbers $i,j$ let $\PNC(i,j)$ denote
the space of polynomials of degree less than or equal to $j$ in $i$
non-commuting variable.
\endproclaim

It should be noted that $I$ is assumed to be in $\PNC(i,j)$ and that
not all $i$ variables need to be represented in every element of
$\PNC(i,j)$.
\medskip

In the rest of this section we will consider a fixed linear subset
$\Cal L$ of $\Cal B$ such that $I\in\Cal L$ and $\Cal L$ generates
$\Cal B$ algebraically. We will let $d$ denote the linear dimension
of $\Cal L$ and consider a fixed basis $\{\ell_1=I$,
$\ell_2,\cdots,\ell_d\}$ for $\Cal L$. Our aim is to get an estimate
of the costs involved in order to use Lemma 2.1. For $k\in\Bbb N$ we
will let $\Cal L^k=\span\{\ell_{i_1}\cdots\ell_{i_k}\mid
i_1,\cdots,i_k\in\{1,\cdots,d\}\}$. Since $I\in\Cal L$ we have $\Cal
L^{k+1}\supseteq\Cal L^k$ so either $\Cal L^k=\Cal L^{k+1}$ or
$\dim(\Cal L^{k+1})\geq\dim(\Cal L^k)+1$. If $\Cal L^k=\Cal L^{k+1}$
then $\Cal L^k=\Cal B$. Hence we get $\Cal B=\Cal L^k$ at least for
$k=\dim(\Cal B)-d+1\,$.
On the other hand this is a very rough estimate relating the linear dimension
of $\Cal L$ to that of $\Cal B$. In fact in order to apply Lemma 2.1 we do
only need to get to the least $k$ such that $ ({\Cal L }^k+\Cal J) =
\Cal B$, hence we define:

\proclaim{2.3 Definition} The semi-simple defect $\sd(\Cal L)$ of $\Cal L$ is
the smallest natural number t such that $$ ({\Cal L}^{t+1}+\Cal J) = \Cal B$$.
\endproclaim

The reason why we think it is worth introducing this term is based on an
example where the matrices in $\Cal L$ generate the upper triangular
matrices. In this case $\sd(\Cal L) \leq (d-1)n$ which for a small $d$ and a
large $n$ is far from $\dim(\Cal B)-d= n(n+1)/2-d\,$.  Using these
definitions and the remarks made we have the following immediate application of
Lemma 2.1.

\proclaim{2.4 Proposition} Let $b\in\Cal B$ then $b\in\Cal J$ if and
only if for any $p\in\PNC(d-1,\sd(\Cal L)+1):$
$$\tr_n\left(bp(\ell_2,\cdots,\ell_d)\right)=0\,.
$$
The semi-simple defect satisfies $\sd(\Cal L) \leq (\dim(\Cal B)-\dim(\Cal L))\,$. 

\endproclaim

The real content of this proposition is mostly that it calls the attention to
the interplay between traces and algebraic properties.  As we shall see below
a well known result of McCoy follows easily from Proposition 2.4, and in order
to provide basic results for the coming sections we formulate and prove some
results based on this type of arguments.

\proclaim{2.5 Proposition} Let $x_1,\cdots,x_s$ be linearly
independent matrices in $M_n(\Bbb C)$ and let
$t=\sd(\span(I,x_1,\cdots,x_s))$, then $t\leq n^2-s$ and the set
$\{x_1,\cdots,x_s\}$ can be triangularized simultaneously if and only
if for any pair $i,j\in\{1,\cdots,s\}$ and any $p\in\PNC(s\,,\;t+1)$
$$\tr_n\left((x_ix_j-x_jx_i)p(x_1,\cdots,x_s)\right)=0\,. \tag$\ast$ $$
\endproclaim

\demo{Proof} The condition is clearly satisfied for sets of
upper triangular matrices.

Next suppose ($\ast$) is valid and let $\Cal B$ denote the algebra
generated by $I=I_{M_n(\Bbb C)}$ and the set $\{x_1,\cdots,x_s\}$, and
let the decomposition from (1) be
 $$\Cal B=M_{n_1}(\Bbb C)\oplus\cdots\oplus M_{n_s}(\Bbb C)\oplus\Cal J\,.
 $$
 
By Proposition 2.4 we find that $x_ix_j-x_jx_i\in\Cal
J$ so the semi-simple part $\Cal C$ of $\Cal B$ is commutative and
$n_1=n_2=\cdots=n_s=1$. This means that for some suitable basis for
$\Bbb C^n;\Cal B$ can be represented as a sub-algebra of the upper
triangular matrices. The reason being that $\Cal B$, as a Lie algebra
under the usual commutator product, must be solvable since $[\Cal
B,\Cal B]\subseteq\Cal J$ and $\Cal J^n=0$. An application of Lie's
Theorem [7, Cor\. A p\. 17] tells that $\Cal B$ is triangularizable.
\hfill$\square$
 \enddemo

The proposition has an immediate corollary

\proclaim{2.6 Corollary} Let $t = \sd( I, x_1,\cdots, x_s )$
then $\{x_1,\cdots,x_s\}$ are si\-mul\-taneous\-ly
tri\-an\-gu\-lar\-iz\-able if for any monomial $x_{i_1}x_{i_2}\cdots
x_{i_m}$ with $m\leq t+3$ and any permutation $\sigma\in\sum_m$
 $$\tr_n(x_{i_1}x_{i_2}\cdots x_{i_m})=
\tr_n\left(x_{i_{\sigma(1)}}x_{i_{\sigma(2)}}\cdots
x_{i_{\sigma(m)}}\right)\,.
 $$
\endproclaim

Another immediate corollary is the following 
\proclaim{2.7 Corollary {\rm[14]}} Let $x,y$ be in $M_n(\Bbb C)$ then
$x$ and $y$ are simultaneously triangularizable if and only if for
any polynomial $p$ in $2$ non-commuting variables $p(x,y)(xy-yx)$ is
nilpotent.
\endproclaim

\demo{Proof} The trace of a nilpotent element vanishes.
\hfill$\square$
 \enddemo

In the rest of this section we will do some computations for 2
matrices in $M_2(\Bbb C)$ and $M_3(\Bbb C)$ respectively. We get some
estimates on how large the semi-simple defect can be in these low
dimensional cases.

\proclaim{2.8 Corollary} Let $x,y\in M_2(\Bbb C)$ then $x$ and $y$
are simultaneously triangularizable if and only if
 $$\tr_2(x^2y^2)=\tr_2\left((xy)^2\right)\,.
 $$
 \endproclaim

\demo{Proof} If $\dim(\span(I,x,y))\leq2$ then $x$ and $y$ commute
and hence they are simultaneously triangularizable. If
$\dim(\span(I,x,y))=3$ then $\Alg(I,x,y)=\{p(x,y)\mid p\in\PNC(2,2)\}$
so the condition in Proposition 2.5 becomes
 $$0=\tr_2\left((xy-yx)p(x,y)\right)\,,  p\in\PNC(2,2)\,.
 $$
 
The only relations which do not vanish automatically are those
where $p$ contains $xy$ or $yx$. In both cases we then get the
condition $\tr_2((xy)^2)=\tr_2(x^2y^2)$.

\hfill$\square$
 \enddemo

The content of Corollary 2.8 is closely related to the description
given by Friedland in [4]. Here it is  proven that $x,y$ in
$M_2(\Bbb C)$ are simultaneously triangularizable if and only if 
 $$
\left(2\tr_2(x^2)-\tr_2(x)^2\right)
\left(2\tr_2(y^2)-\tr_2(y)^2\right)=\left(2\tr_2(xy)-\tr_2(x)\tr_2(y)\right)^2\,.
$$

If one applies Cayley-Hamiltons Theorem to $x^2$, $y^2$ and $(xy)^2$
in the relation presented in Corollary 2.8, one can obtain Friedlands
relation quite easily.
\medskip

We will now investigate the case where $x, y$ are in $M_3(\Bbb C)$.  An
immediate application of Proposition 2.5 in order to determine whether $x, y$
are simultaneously triagularizable or not would involve traces of polynomials
in $x, y$ of total degree $9$, $(9=(3^2-3)+3)\,.$ It should however be
remarked that the problems we are facing when we have to compute traces of
polynomials are reduced considerably by the fact that the trace is invariant
under cyclic permutations of monomials. We will take a closer look into this
problem for $n = 3$ in order to see that degree $6$ suffices whereas $5$ does
not.  We start with an example demonstrating that polynomials of degree $5$
are not sufficient to settle the question for $x, y$ in $M_3(\Bbb C)$.  
\medskip

\proclaim{2.9 Example}  
 
\endproclaim
Let $x,y$ in $M_3(\Bbb C)$ be given by
 $$x=\pmatrix
0&\;\;0&\;\;0\\
0&\;\;2&\;\;0\\
0&\;\;0&\;\;1+i\sqrt3\endpmatrix
\qquad y=
\pmatrix
0&\;\;1&\;\;0\\
0&\;\;0&\;\;1\\
1&\;\;0&\;\;0
\endpmatrix\,.
$$
 Clearly $\Alg(x,y)=M_3(\Bbb C)$ so the matrices can not be
triangularized simultaneously. On the other hand we have
 $$\forall p\in\PNC(2,3):\;\tr_3\left((xy-yx)p(x,y)\right)=0\,.\tag2
 $$ This is easily seen once it is observed that for any monomial
$m=x^{i_1}y^{j_1}\cdots x^{i_k}y^{j_k}$ we have $\tr_3(m)=0$ unless
$j_1+\cdots+j_k\in3\Bbb Z$. 

Hence we only have to evaluate 
 $$\tr_3\left((xy-yx)y^2\right),\;\,
\tr_3\left((xy-yx)xy^2\right),\;\,
\tr_3\left((xy-yx)yxy\right),\;\, \tr_3\left((xy-yx)y^2x\right)\,.
 $$ The first of these terms vanishes because $\tr_3(ab)=\tr_3(ba)$. Using
this property over and over again, the other expressions are seen to vanish
because they can be reduced to
 $$
\tr_3(y^2xyx-y^3x^2)\,,\quad0\,,\quad
\tr_3(y^3x^2-y^2xyx)\,,
$$
and a simple computation shows that $\tr_3(y^3x^2)=\tr_3(y^2xyx)$.
\hfill$\square$
\medskip

\proclaim{2.10 Proposition} Let $x,y$ be in $M_3(\Bbb C)$ then $x$ and
$y$ can be triangularized simultaneously if and only if for any
monomial $m=x^{i_1}y^{j_1}x^{i_2}y^{j_2}x^{i_3}y^{j_3}$ of degree $6$
or less $\tr_3(m)=\tr_3(x^{(i_1+i_2+i_3)}y^{(j_1+j_2+j_3)})$.
\endproclaim

\demo{Proof} Suppose the condition is satisfied then by Corollary 2.6 it
suffices to show that the algebra, say $\Cal B$, generated by $\{I,x,y\}$ is
spanned by polynomials $p(x,y)$ of degree $4$ or less. Let $\Cal B=\Cal
C\oplus\Cal J$ be the algebraic Wedderburn decomposition (1) where $\Cal C$
denotes the semi-simple summand. Then $\Cal C$ can be of one of the forms
$\Bbb C\oplus\Bbb C\oplus\Bbb C$, $\Bbb C\oplus M_2(\Bbb C)$ or $M_3(\Bbb C)$,
so we will have to exclude the latter $2$ under the assumptions made on $x$
and $y$. Let us first suppose $\Cal C=\Bbb C\oplus M_2(\Bbb C)$. Then for
suitable elements $\lambda_x,\lambda_y\in\Bbb C$, $m_x,m_y\in M_2(\Bbb C)$ and
$j_x,j_y$ in $\Cal J$ we have
 $$x=\lambda_x+m_x+j_x\qquad y=\lambda_y+m_y+j_y\,.
 $$
 
Clearly $m_x,m_y$ are not simultaneously triangularizable since
they must ge\-ne\-ra\-te $M_2(\Bbb C)$. By
Corollary 2.8 we have
$\tr_3(m_x^2m_y^2)\neq\tr_3\left((m_xm_y)^2\right)$ and then
 $$\tr_3(x^2y^2)-\tr_3(xyxy)=\tr_3\left(m_x^2m_y^2-(m_xm_y)^2\right)\neq0\,,
 $$
 which contradicts the assumptions made on $x$ and $y$.

We may therefore suppose that $\Cal B=M_3(\Bbb C)$ and  we prove below 
 that under this assumption one has
$$
M_3(\Bbb C)=\Alg(I,x,y)=\{p(x,y)\mid
p\in\PNC(2,4)\}\,.\tag$\ast\ast$
$$ By Corollary 2.6 $x$ and $y$ are then simultaneously triangularizable and
the assumption $\Cal B=M_3(\Bbb C)$ can not be true.

The set $\{I,x,y,x^2,y^2,xy,yx\}$ has 7 elements. If it is linearly independent 
then since   $\Cal B=M_3(\Bbb C)$  the
linear dimension of $\{p(x,y)\mid p\in\PNC(2,4)\}$ must be at least
$9$ and hence $9$ so we have ($\ast\ast$). On the other hand if
the set is not linearly independent we will prove that there exists a
linear relation of the type
 $$\lambda_1+\lambda_2x+\lambda_3y+
\lambda_4x^2+\lambda_5y^2+\lambda_6xy+\lambda_7yx=0\,,
$$
where either $\lambda_6\neq0$ or $\lambda_7\neq0$. Let us postpone
the proof of this and see how this can be used to prove the
proposition. Suppose for instance $\lambda_7\neq0$ and say
$\lambda_7=1$ then
 $$yx=-\lambda_1-\lambda_2x-\lambda_3y-\lambda_4x^2-\lambda_5y^2-\lambda_6xy\,,
 $$
 so any polynomial $p(x,y)$ can be written as a sum of products of
polynomials in just one variable
 $$p(x,y)=q_0(x)+r_0(y)+\sum_{i=1}^tq_i(x)r_i(y)\,.
 $$ By Cayley-Hamiltons Theorem we then get
$p(x,y)=\sum_{i=0}^2\sum_{j=0}^2\alpha_{ij}x^iy^j$ and it turns out that $\Cal
B=\{p(x,y)\mid p\in\PNC(2,4)\}$, and the argument is completed. 

  Let us then assume that the set $\{1,x,y,x^2,y^2,xy,yx\}$ is linearly
dependent and satisfies a non trivial linear relation as above with $\lambda_6
= \lambda_7 = 0$. If $\lambda_4 = \lambda_5 = 0$ also then $x$ and $y$
commutes so we may and will assume that $\lambda_5 = 1$. let $z=y-x$ then
$y=x+z$ and the linear relation becomes - when expressed in $x$ and $z$:
 $$\lambda_1+(\lambda_2+\lambda_3)x+\lambda_3z+(1+\lambda_4)x^2+z^2+xz+zx=0\,,
 $$
 and we are back in a situation already covered.

\hfill$\square$
\enddemo
\subhead{3. Property L and tensor products}\endsubhead

Let $x,y$ be $2$ matrices in $M_n(\Bbb C)$, then according to [12,
15] $x$ and $y$ are said to have Property $L$ if there exist sets
$(s_1,\cdots,s_n)$ $(t_1,\cdots,t_n)$ of complex numbers such that
for $\lambda,\mu\in\Bbb C$ the characteristic roots of $(\lambda
x+\mu y)$ is the set $\{\lambda s_i+\mu t_i\mid 1\leq i\leq n\}$.
Notably O\. Taussky [16, 17] and T\. Laffey [10,11] have obtained
results on such pairs of matrices. Especially sufficient conditions
which together with Property $L$ implies that $x$ and $y$ are
simultaneously triangularizable have been searched for. The following
example comes from Wielandt [19] but was used in [15] to show that
there exist 2 nilpotent matrices $x$ and $y$ in $M_3(\Bbb C)$ such
that every element $\lambda x+\mu y$ in the matrix pencil 
is nilpotent and $\alg(x,y)=M_3(\Bbb C)$. The
matrices $x$ and $y$ then do have Property $L$, but they are not
simultaneously triangularizable.
\medskip

\proclaim {3.1  Example, \rm{[15, 19]}}  

\endproclaim
Let $x=\left(\smallmatrix
0&\;0&\;0\\
1&\;0&\;0\\
0&\;1&\;0\endsmallmatrix\right)$ and 
$y=\left(\smallmatrix
0&\;1&\;\phantom{-}0\\
0&\;0&\;-1\\
0&\;0&\;\phantom{-}0\endsmallmatrix\right)$ then $\lambda x+\mu y$ is nilpotent for
$\lambda,\mu\in\Bbb C$, and $\alg(x,y)=M_3(\Bbb C)$.

\proclaim {3.2 Definition}
 Let $\Cal S\subseteq M_n(\Bbb C)$ be a set of
matrices such that for each $a\in\Cal S$  the characteristic roots of
$a$ are equipped with a given numbering
$(\lambda_1^a,\cdots,\lambda_n^a)$ and let $k\in\Bbb N$. Then $\Cal S$ 
is said to have the property $kL$ if
for any set $(a_1,\cdots,a_j)$ from $\Cal S$ and any set
$(x_1,\cdots,x_j)$ from $M_k(\Bbb C)$ and $t\in\Bbb C$
$$\aligned
&\det\left(t(I_k\otimes I_n)-(x_1\otimes a_1+\cdots+x_j\otimes
a_j)\right)\\
&\qquad\qquad\qquad =\underset{i=1}\to{\overset n\to\prod}\det\left(
tI_k-(\lambda_i^{a_1}x_1+\cdots+\lambda_i^{a_j}x_j)\right)\,.
\endaligned\tag$\ast$
$$
\endproclaim
It is clear from the definition that the characteristic roots of the
sum $x_1\otimes a_1+\cdots+x_j\otimes a_j$ are grouped as the
disjoint union of the sets 
 $$
\text{ characteristic roots of }\; (\lambda_i^{a_1}x_1+\cdots+\lambda_i^{a_j}x_j)\,, 1 \leq i\leq n
 $$
 when $\Cal S$ has Property $kL$. We have chosen
\lq\lq$kL$\rq\rq\,rather than \lq\lq$Lk$\rq\rq\;because \lq\lq
Property $Lk$\rq\rq\;has been used in [18] to mean the generalization
of \lq\lq Property $L$\rq\rq\;to sets of $k$ matrices.

The shift from \lq\lq Property $L$\rq\rq\;to \lq\lq Property
$kL$\rq\rq\;corresponds to a shift in coefficients from scalars to
matrices in matrix pencils
 $$\lambda a+\mu b\to x\otimes a+y\otimes b\,.
 $$
 
This method has been very fruitful in $K$-theory and in many so-called
non-commutative theories. Further we think it fits well with some of the early
works in matrix theory by Kronecker [5,12].  After all the tensor product was
known and used as the Kronecker product in matrix theory long before its
general algebraic nature was described.  The following lemma needs no written
proof, but the observation has to be made.

\proclaim{3.3 Lemma} 
If a set $\Cal S$ in $M_n(\Bbb C)$ has Property
$kL$ then it has Property $(k-1)L$.
\endproclaim

The following proposition shows how the logarithm can give an
alternative characterization of Property $kL$ via the trace rather
than the determinant.

\proclaim{3.4 Proposition} Let $\Cal S\subseteq M_n(\Bbb C)$ be such
that for each $a\in\Cal S$ there is a given numbering
$(\lambda_1^a,\cdots,\lambda_n^a)$ of its characteristic roots. Then
$\Cal S$ has property $kL$ if and only if for any $j$ in $\Bbb N$ and
any sets $(a_1,\cdots,a_j)$ from $\Cal S$, $(x_1,\cdots,x_j)$
from $M_k(\Bbb C)$ and $m\in\Bbb N$, $1\leq m\leq nk$:
$$
\tr_{nk}\left((x_1\otimes a_1+\cdots+x_j
\otimes a_j)^m\right)=\sum_{i=1}^n
\tr_k\left((\lambda_i^{a_1}x_1+\cdots+\lambda_i^{a_j}x_j)^m\right)\tag$\ast$
$$

If $\Cal S$ has property $kL$ then ($\ast$) is true for all $m\in\Bbb
N$.
\endproclaim

\demo{Proof} Suppose $\Cal S$ has property $kL$, then from Definition
3.2 we get for $z\in\Bbb C$
 $$\aligned
&\det\left(I_k\otimes I_n-z(x_1\otimes a_1+\dots+x_j\otimes
a_j)\right)\\
 &\qquad\qquad\qquad =
\underset{i=1}\to{\overset n\to\prod}
\det\left(I_k-z(\lambda_i^{a_1}x_1+
\cdots+\lambda_i^{a_j}x_j)\right)\,.\endaligned\tag$\ast\ast$
 $$ For $z$ near $0$, the main branch of the logarithm can be applied on both
sides of ($\ast\ast$) and we get since $\Log(\det(y))=\tr_r(\Log(y))$ when $y$
is near $I$ in any $M_r(\Bbb C)$
 $$\align
&\tr_{kn}\left(\Log(I_k\otimes I_n-z(x_1\otimes a_1+\cdots+x_j\otimes
a_j))\right)\\
&\qquad\qquad\qquad  =
\sum_{i=1}^n\tr_k\left(\Log\left(I_k-z
(\lambda_i^{a_1}x_1+\cdots+\lambda_i^{a_j}x_j)\right)\right)\,.\endalign
 $$ Both sides of this equation can be expressed as power series. A comparison
of terms yields that ($\ast$) holds for all $m\in\Bbb N$.  \smallskip

Now suppose ($\ast$) is valid for all natural numbers $m$ less than or equal
to $nk$. Then choose diagonal matrices $d_1,\cdots,d_j$ in $M_n(\Bbb C)$ such
that
$$d_\ell=\Delta(\lambda_1^{a_\ell},\lambda_2^{a_\ell},\cdots,\lambda_n^{a_\ell})\,.$$
Clearly the set $(d_1,\cdots,d_j)$ has Property $kL$ so for any $m$ in
$\Bbb N$
 $$\tr_{nk}\left((x_1\otimes d_1+
\cdots+x_j\otimes d_j)^m\right)=
\sum_{i=1}^n\tr_k\left((\lambda_i^{a_1}
x_1+\cdots+\lambda_i^{a_j}x_j)^m\right)\,.
$$ Hence in order to prove that ($\ast$) holds for all $m$ we just have to
note that by assumption and the computations just made:
 $$\align
&\text{for }\; 1\leq m\leq nk:\quad \tr_{nk}\left((x_1\otimes
a_1+\cdots+x_j\otimes a_j)^m\right)\\
&\qquad\qquad\qquad=\tr_{nk}\left((x_1\otimes d_1+\cdots+x_j\otimes
d_j)^m\right)\,.\endalign
 $$ This means by Newtons formulae that the characteristic polynomia for
 $x_1\otimes a_1+\cdots+x_j\otimes a_j$ and $x_1\otimes d_1+\cdots+x_j\otimes
 d_j$ are identical and hence that for any $m$ in $\Bbb N$
 $$\align
\tr_{nk}\left((x_1\otimes a_1+\cdots+x_j\otimes a_j)^m\right)&=
\tr_{nk}\left((x_1\otimes d_1+\cdots+x_j\otimes d_j)^m\right)\\
&=\sum_{i=1}^n\tr_k
\left((\lambda_i^{a_1}x_1+\cdots+\lambda_i^{a_j}x_j)^m\right)\,.\endalign
$$

On the other hand if all equations of the type ($\ast$) are true for
all $m\in\Bbb N$ it follows that ($\ast\ast$) holds in some ball
around zero and Property $kL$ is established.\hfill$\square$  
\enddemo

The following lemma is included in order to
prepare an application of results from Section 2. 

\proclaim{3.5 Lemma} If a set $\Cal S\subseteq M_n(\Bbb C)$ has a
property $kL$ then so does $\Cal S\cup I_{M_n(\Bbb C)}$.
\endproclaim

\demo{Proof} Let $a_1,\cdots,a_j$ be in $\Cal S$,
$x_0,x_1,\cdots,x_j$ be in $M_k(\Bbb C)$ and $z\in\Bbb C$ then for
$z\in B(0,r)$ for some $r>0$ we have $(I_k-zx_0)$ is invertible in
$M_k(\Bbb C)$ so 
$$\align
&\det\left(I_k\otimes I_n-(z x_0\otimes I_n+x_1\otimes
a_1+\cdots+x_j\otimes a_j)\right)\\
=&\det\left((I_k-zx_0)\otimes I_n\right)\det\left(I_k\otimes
I_n - \left((I_k-zx_0)^{-1}x_1\otimes
a_1+\cdots\right.\right.\\
&\left.\left.\qquad\qquad\qquad\qquad
\qquad\qquad\qquad\qquad\qquad\qquad\qquad
+(I_k-zx_0)^{-1}x_j\otimes a_j\right)\right)\\ 
=&\det(I_k-zx_0)^n\underset{i=1}\to{\overset
n\to\prod}\det
\left(I_k-(I_k-zx_0)^{-1}
(\lambda_i^{a_1}x_1+\cdots+\lambda_i^{a_j}x_j)\right)\\
=&\underset{i=1}\to{\overset
n\to\prod}
\det\left(I_k-(zx_0+\lambda_i^{a_1}x_1+
\cdots\lambda_i^{a_j}x_j)\right)\,.
\endalign
 $$ Hence this identity is true for $z=1$ as well and the lemma
follows.\hfill$\square$ 
\enddemo

We can now state the main result of this section 

\proclaim{3.6 Theorem} Let $\Cal S\subseteq M_n(\Bbb C)$ be a set and
$k=\sd(\span(\Cal S \cup I))+3$.
The matrices in  $\Cal S$ are  simultaneously triangularizable if and only if
$\Cal S$ has Property $kL$. 
Let $s=\dim(\span(\Cal S))\,$ then $k \leq n^2-s+3$.

\endproclaim

\demo{Proof} It is obvious that Property $kL$ is necessary, so let us
assume that $\Cal S$ has Property $kL$ then by Lemma 3.5 $\{\Cal S\cup I\}$
has Property $kL$.  According to Corollary 2.6 it suffices
to prove that for any monomial $a_1\cdots a_m$ of degree $k$ or less of
elements from $\Cal S$ and any permutation $\sigma\in\sum_m$ we have
 $$\tr_n\left(a_1\cdots a_m\right)=\tr_n\left(a_{\sigma(1)}\cdots
a_{\sigma(m)}\right)\,.
 $$ Since $I\in\{\Cal S\cup I\} $ we can replace this statement with the
following.  The set $\Cal S$ is simultaneously triangularizable if for any set
$a_1,\cdots,a_k$ from $\{\Cal S \cup I\} $ and any $\sigma$ in $\sum_k$
$$
\tr_n(a_1,\cdots a_k)=\tr_n( a_{\sigma(1)}\cdots a_{\sigma(k)})\,.
$$

In order to see that Property $kL$ implies the statement above we let
$(e_{ij})$ denote the matrix-units in $M_k(\Bbb C)$, $\sigma$ a
permutation in $\sum_k$ and let $a_1,\cdots,a_k$ be elements from $\{\Cal
S\cup I \}$. Then we define 
$$\align
u&=e_{12}\otimes a_1+e_{23}\otimes a_2+\cdots+e_{\left(k-1 \right)k}\otimes
a_{k-1}+e_{k1}\otimes a_k\\
u^\sigma&=e_{12}\otimes a_{\sigma(1)}+e_{23}\otimes
a_{\sigma(2)}+\cdots +e_{k1}\otimes a_{\sigma(k)}\,.\endalign 
 $$ A computation -- which uses that the trace is invariant under cyclic
permutations of products -- shows that
 $$\align
\tr_{nk}(u^k)&=k\,\tr_n(a_1a_2\cdots a_k)\\
\tr_{nk}\left((u^\sigma)^k\right)&=k\,\tr_n(a_{\sigma(1)}a_{\sigma(2)}\cdots
a_{\sigma(k)})\,.\endalign
 $$ Analogously we define for $1\leq i\leq n,v_i.v_i^\sigma$ in $M_k(\Bbb C)$
by
 $$\align
v_i&=\lambda_i^{a_1}e_{12}+
\lambda_i^{a_2}e_{23}+\cdots+\lambda_i^{a_k}e_{k1}\\
v_i^\sigma&=\lambda_i^{a_{\sigma(1)}}e_{12}+
\lambda_i^{a_{\sigma(2)}}e_{23}+\cdots+
\lambda_i^{a_{\sigma(k)}}e_{k1}\endalign
$$
and we find as above that
 $$\tr_k(v_i^k)=k\,\lambda_i^{a_1}\cdots
\lambda_i^{a_k}=\tr_k\left((v_i^\sigma)^k\right)\,.
$$ Since $\{\Cal S \cup I \}$ has Property $kL$ we then get
$$\align
\tr_n(a_1a_2\cdots a_k)&=\frac1k\tr_{nk}(u^k)\\
&=\frac1k\sum_{i=1}^n\tr_k(v_i^k)\\
&=\frac1k\sum_{i=1}^n\tr_k\left((v_i^\sigma)^k\right)\\
&=\frac1k\tr_{nk}\left((u^\sigma)^k\right)\\
&=\tr_n(a_{\sigma(1)}a_{\sigma(2)}\cdots a_{\sigma(k)})\,,\endalign
$$
 and the theorem follows.\hfill$\square$  
\enddemo

\subhead{4. Invertibility preserving mappings}\endsubhead

In an earlier article [3] we have obtained a result on continuous
linear mappings from a Banach algebra into $M_n(\Bbb C)$ which
preserves invertibility. The result is a generalization of  the
Gleason-Kahane-\`Zelazko Theorem [6,8]. In the situation here, we can quote
[3] in the following way.

\proclaim{4.1 Theorem {\rm ([3])}} Let $\Cal A\subseteq M_k(\Bbb C)$
be a unital algebra and $\varphi:\Cal A\to M_n(\Bbb C)$ a linear
mapping satisfying $\varphi(I_\Cal A)=I_{M_n(\Bbb C)}$. Let $\Cal
A_{\text{inv}}$ denote the set of invertibles in $\Cal A$ then
$\varphi(\Cal A_{\text{inv}})\subseteq \GL_n(\Bbb C)$ if and only if
for any $a$ in $\Cal A$ and any $m$ in $\Bbb N$
 $$\tr_n\left(\varphi(a^m)\right)=\tr_n\left(\varphi(a)^m\right)\,.
 $$
 \endproclaim

There is an immediate corollary which is quite useful.

\proclaim{4.2 Corollary {\rm [(3])}} If $\varphi(\Cal
A_{\text{inv}})\subseteq\GL_n(\Bbb C)$ then for any $a,b$ in $\Cal A$
and any $k$ in $\Bbb N$
\roster
 \item"(i)"
$\tr_n(\varphi(ab))=\tr_n(\varphi(a)\varphi(b))=\tr_n(\varphi(ba))$.
\smallskip
 \item"(ii)"  $\tr_n(\varphi(a)^k\varphi(b))
=\tr_n(\varphi(a^kb))=\tr_n(\varphi(a^k)\varphi(b))$
\smallskip
 \item"(iii)"  $\det(\varphi(a)\varphi(b))=\det(\varphi(ab))$.
 \endroster
\endproclaim

\demo{Proof} The relation (ii) is not stated explicitly 
 in [3], but follows
easily from Theorem 4.1. Let $z \in \Bbb{C}
$ and $x=(a+zb)$ then by Theorem 4.1
 $$\tr_n\left(\varphi(x^{k+1})\right)=\tr_n\left(\varphi(x)^{k+1}\right)\,.
 $$ The relation (i) and the trace properties for $\tr_n$ shows that for the
coefficient to $z$ in the identity right above we get
 $$k\tr_n\left(\varphi(a^kb)\right)
=k\tr_n\left(\varphi(a)^k\varphi(b)\right)\,.
$$ An application of (i) once more gives (ii).\hfill$\square$
\enddemo

One of the main motivations to look into the problem of trying to
describe invertibility preserving linear mappings was the set of
lecture notes [9] where Kaplansky addresses this problem. As
mentioned in the introduction Kaplansky had the impression that the
-- at the time of the notes -- quite recent results [6,8] by Gleason
and Kahane \& \`Zelazko might be generalized. Of the various articles
we have found in this area of research especially the work by Aupetit
[1] has been fruitful to us. The content in this article as well
as the one in [3] is very much influenced by [1]. As mentioned in the
introduction Kaplansky was a bit too optimistic in hoping that
invertibility preserving mappings should be Jordan homomorphisms. The
following examples demonstrate what sort of obstacles we have found.
\medskip

\proclaim {4.3 Examples.}
\endproclaim

\flushpar{\bf A: \quad} Let $\Cal A$ be the diagonal matrices in
$M_3(\Bbb C)$, $\Cal T$ be the upper triangular matrices in $M_3(\Bbb
C)$ and let $\varphi:\Cal A\to\Cal T$ be given by
 $$\varphi(\Delta(a,b,c))=\pmatrix
a&\quad 0&\quad0\\
0&\quad b&\quad a-b\\
0&\quad 0&\quad c\endpmatrix\,.
 $$
 
This $\varphi$ is unital and invertibility preserving but
$\varphi(\Delta(1,0,0)^2)=\varphi(\Delta(1,0,0))\neq\varphi(\Delta(1,0,0))^2$. On
the other hand $\varphi$ is a homomorphism modulo the radical.
\medskip

\flushpar {\bf B: \quad} Let $\Cal A$ be the diagonal matrices in
$M_3(\Bbb C)$ and $\varphi:\Cal A\to M_3(\Bbb C)$ be given by
 $$\varphi\left(\Delta(a,b,c)\right)=\pmatrix
a&\quad (c-a)&\quad 0\\
(b-a)&\quad a&\quad (a-c)\\
0&\quad (b-a)&\quad a\endpmatrix\,,
 $$
 then $\det\varphi(\Delta(a,b,c))=a^3$ and $\varphi(I)=I$, but
$\varphi$ is not a homomorphism and by Example 3.1 the
algebra generated by $\varphi(\Cal A)$ equals $M_3(\Bbb C)$ which is
semi-simple and then has no radical. Hence not even modulo the radical
do we get a homomorphism.
\medskip

\flushpar{\bf C: \quad} Let $\Cal A=M_2(\Bbb C)$ and $\varphi:\Cal A\to
M_6(\Bbb C)$ be given by
 $$\varphi\pmatrix
a&\;b\\
c&\;d\endpmatrix=
\pmatrix
aI_3&\quad bI_3\\
&\\
cI_3&\quad
\pmatrix
d&\;b&\;0\\
a-d&\;d&\;-b\\
0&\;a-d&\;d\endpmatrix\endpmatrix\;.
 $$
 
Clearly by Example 3.1 for $x\in M_2(\Bbb C)$
$\det(\varphi(x))=(\det(x))^3$ and  $\varphi(I_2)=I_6$. The algebra
generated by $\varphi(\Cal A)$ equals $M_6(\Bbb C)$, so the Jacobson
radical vanishes and $\varphi$ can not be a Jordan homomorphism
modulo the radical.
\medskip

We will now turn to some positive results. The first is closely
related to Aupetit's result [1].

\proclaim{4.4 Theorem} Let $\Cal A\subseteq M_h(\Bbb C)$ be a
unital algebra and $\varphi:\Cal A\to M_n(\Bbb C)$ a unital
invertibility preserving mapping. If $\varphi(\Cal A)$ is an algebra
then $\varphi$ is a Jordan homomorphism modulo the Jacobson radical.
\endproclaim

\demo{Proof} Let $x$ be in $\Cal B=\alg\varphi(\Cal A)=\varphi(\Cal
A)$ and $a$ in $\Cal A$ then for $k\in\Bbb N$ we get from Corollary
4.2 (ii), and the assumptions made that there exists $b$ in $\Cal A$
such that $x=\varphi(b)$ and hence 
 $$\tr_n\left((\varphi(a)^k-\varphi(a^k))x\right)=
\tr_n\left(\varphi((a^k-a^k)b)\right)=0\,.
 $$ By Lemma 2.1 $\varphi(a)^k-\varphi(a^k)$ belongs to the Jacobson radical
so $\varphi$ is a Jordan homomorphism modulo the Jacobson
radical.\hfill$\square$
\enddemo

Analogously to the results in section 3 we can obtain sufficient
conditions if we demand that $\varphi\otimes\text{ id}$ on $\Cal
A\otimes M_k(\Bbb C)\to M_{nk}(\Bbb C)$ is invertibility preserving,
hence we define.

\proclaim{4.5 Definition} Let $\Cal A$ be a unital algebra in
$M_h(\Bbb C)$ and $\varphi:\Cal A\to M_n(\Bbb C)$ a unital
invertibility preserving linear mapping. For $k$ in $\Bbb N$,
$\varphi$ is said to be $k$-invertibility preserving if
$\varphi\otimes\text{ id}_k:\Cal A\otimes M_k\to M_n\otimes M_k$
preserves invertibility of elements.
\endproclaim

\proclaim{4.6 Theorem} Let $\Cal A$ be a unital algebra in $M_h(\Bbb C)$,
$\varphi:\Cal A\to M_n(\Bbb C)$ a unital linear mapping, $\Cal B$ the algebra
generated by $\varphi(\Cal A)$, $t=\sd(\varphi(\Cal A))$ and $k=t+3$. The
mapping $\varphi$ is $k$-invertibility preserving if and only if $\varphi$ is
a homomorphism modulo the Jacobson radical. The semi-simple defect $t$ of
$\varphi(\Cal A)$ satisfies $t \leq \dim(\Cal B) - \dim(\varphi(\Cal A))$.
\endproclaim

\demo{Proof} If $\varphi$ is a homomorphism modulo the radical then so is
$\varphi\otimes\text{ id}_m$ for any natural number $m$, and $\varphi$ is
$m$-invertibility preserving for all $m$.  
\smallskip

Let us now assume that $\varphi$ is $k$-invertibility preserving and let $u,v$
be in $\Cal A$. In order to prove that $x=(\varphi(uv)-\varphi(u)\varphi(v))$
belongs to Jacobson radical $J$ of $\Cal B$ we have by Proposition 2.3 to
prove that for any $m\leq k-2$ and any set $a_1,\cdots,a_m$ from $\Cal A$ we
have
 $$\tr_n(x\varphi(a_1)\cdots\varphi(a_m))=0\,.
 $$ Since $I_n=\varphi(I_{\Cal A})$ it is sufficient to prove that
$\tr_n(x\varphi(a_1)\cdots\varphi(a_{k-2}))=0\,$ for all sets
$a_1,\cdots,a_{k-2}$ from $\Cal A$.

We start by proving that for any set $a_1,\cdots,a_k$ from $\Cal A$
we have
 $$\tr\left(\varphi(a_1)\cdots\varphi(a_k)\right)=\tr\left(\varphi(a_1\cdots
a_k)\right)
 $$
 in order for do so we let $(e_{ij})$ denote a set of matrix units in
$M_k(\Bbb C)$ and define
 $$u=a_1\otimes e_{12}+a_2\otimes e_{23}+\cdots+a_k\otimes e_{k1}\in\Cal
A\otimes M_k(\Bbb C)\,.
 $$ In analogy with the computations made in the proof of Theorem 3.6 we get
via Corollary 4.2 (i) applied to $\varphi$
$$ \tr_n\left(\varphi(a_1\cdots a_k)\right)=
\frac1k\tr_{nk}\left((\varphi\otimes\text{ id}_k)(u^k)\right)$$ then since
$\varphi$ is $k$ invertibility preserving we get from Theorem 4.1 that
 $$ \frac1k\tr_{nk}\left((\varphi\otimes\text{
id}_k)(u^k)\right)=\frac1k\tr_{nk}\left(((\varphi\otimes\text{
id}_k)(u))^k\right) $$ 
but the arguments from the proof of Theorem 3.6 applies again and
$$\frac1k\tr_{nk}\left(((\varphi\otimes\text{ id}_k)(u))^k\right) =
\tr_n\left(\varphi(a_1)\cdots \varphi(a_k)\right)\,.$$
Hence for any set $a_1,\cdots,a_{k-2}$ we get
 $$\align
&\tr_n\left(x\varphi(a_1)\varphi(a_2)\cdots\varphi(a_{k-2})\right)\\
=&\tr_n\left((\varphi(uv)\varphi(I)-\varphi(u)\varphi(v))
\varphi(a_2)\cdot\varphi(a_3)\cdots\varphi(a_{k-2})\right)\\
=&\tr_n\left(\varphi((uvI-uv)a_2\cdots a_{k-2}\right)\\
=&0\,,\endalign
 $$
 and the theorem follows.\hfill$\square$ \enddemo 
\medskip
\proclaim {4.7 Remark} 
\endproclaim 
It is well known that transposition on $M_n(\Bbb C)$ is an anti-automorphism,
hence it follows from the theorem that transposition can not be
$k$-invertibility preserving for $k \geq 3$.  On the other hand the following
example for $n = 2$ can be used for all $n \geq 2$ to show that transposition
on $M_n(\Bbb C)$ is never $2$-invertibility preserving.  Let $x \in M_2(\Bbb
C) \otimes M_2(\Bbb C)$ be given by $ x=e_{11}\otimes e_{11}+ e_{21}\otimes
e_{12} + e_{12}\otimes e_{21} + e_{22}\otimes e_{22}$ then x is invertible and
for $\varphi$ as transposition we get $\varphi_2(x)=e_{11}\otimes e_{11}+
e_{12}\otimes e_{12} + e_{21}\otimes e_{21} + e_{22}\otimes e_{22}$ which is
not invertible.

In the proof above we can manage products of length $k$ made by elements from
$\varphi(\Cal A)$. It is expected that compensation in the size of products
manageable when involving $k\times k$ matrices ought to grow like $k^2$ rather
than linearly in $k$. For $k = 2$ it is possible to do much better as the
following proposition shows. We know that Proposition 4.8 has a
generalization to matrices of arbitrary size, but we have not been able to
find a nice description of a general result.
\medskip

\proclaim{4.8 Proposition} Let $\Cal A\subseteq M_h(\Bbb C)$ and
$\varphi:\Cal A\to M_n(\Bbb C)$ be a unital linear map. If $\varphi$
is $2$-invertibility preserving then for all
$\{a,b,c,d\}\subseteq\Cal A$ and all $i,j\in\Bbb N$
\roster
 \item"(i)"
$\tr_n(\varphi(ab^icd^j))=\tr_n(\varphi(a)\varphi(b)^i\varphi(c)\varphi(d)^j)$
\smallskip
 \item"(ii)" $\tr_n(\varphi(ab)^i)=\tr_n((\varphi(a)\varphi(b))^i)$. 
 \endroster
\endproclaim

\demo{Proof} For $s,t$  in $\Bbb C$, $u$ in $\Bbb C\setminus\{0\}$
and $a,b,c,d$ in $\Cal A$ we define $e=u^{-1}(I-td)$ in $\Cal A$ and
$x$ in $\Cal A\otimes M_2(\Bbb C)$ by 
 $$x=\pmatrix
1-sb&\quad c\\
a&\quad e\endpmatrix\,.
 $$
 
Let  $s$ be chosen in a ball $B(0,r)$ such that $(1-sb)$ is
invertible for every s from this ball. We can then define a $2$
invertibility preserving unital map $\psi:\Cal A\to M_n(\Bbb C)$ by
$$ y\in\Cal A:\quad \psi(y)=\varphi\left(y(1-sb)\right)\varphi(1-sb)^{-1}\,.$$

The following identity is straight forward and verifies that $\psi$ is
2-invertibility preserving since for 
$\psi_2=\psi\otimes\text{ id}_{M_2(\Bbb
C)}$ we get 
$$\det(\varphi_2(x))=\det\left(I-s\varphi(b)\right)^2\det(\psi_2(
\pmatrix
1&\quad c(1-sb)^{-1}\\
a(1-sb)^{-1}&\quad e(1-sb)^{-1}\endpmatrix))\,.
$$
 
Let $w\in\Cal A\otimes M_2(\Bbb C)$ be given by
 $$w=\pmatrix
1&\quad 0\\
-a(1-sb)^{-1}&\quad 1\endpmatrix\,,
 $$
 then by Corollary 4.2 (iii) applied to  multiplication by $\psi_2(w)$
from the left gives
 $$\align
\det(\varphi_2(x))&=\det\left(I-s\varphi(b)\right)^2\det(\psi_2(\pmatrix
1&\quad c(1-sb)^{-1}\\
0&\quad (e-a(1-sb)^{-1}c)(1-sb)^{-1}\endpmatrix))\\
&=\det\left(1-s\varphi(b)\right)\det(\varphi\left(e-a(1-sb)^{-1}c\right))\,.
\endalign
 $$ The same type of row operations applied to $\varphi_2(x)$ yields
 $$\det\left(\varphi_2(x)\right)=\det\left(
\pmatrix
1-s\varphi(b)&\quad \varphi(c)\\
0&\quad \varphi(e)-\varphi(a)(I-s\varphi(b))^{-1}\varphi(c)\endpmatrix\right)\,,
 $$
 and therefore
 $$\align
&\det\left(u^{-1}(I-t\varphi(d))-\varphi(a(1-sb)^{-1}c)\right)\\
=&\det\left(u^{-1}(I-t\varphi(d))
-\varphi(a)(I-s\varphi(b))^{-1}\varphi(c)\right)\,.\endalign
$$ So
 $$\align
&\det\left(I-u\varphi(a(1-sb)^{-1}c\right)\left(1-t\varphi(d))^{-1}\right)\\
=&\det\left(I-u\varphi(a)(I-s\varphi(b))^{-1}
\varphi(c)(I-t\varphi(d))^{-1}\right)\,.\endalign
$$ This identity can be extended to $u=0$ and we find by differentiation
and evaluation at $u=0$ that
 $$\align
&\tr_n\left(\varphi(a(1-sb)^{-1}c)(1-t\varphi(d))^{-1}\right)\\
=&\tr_n\left(\varphi(a)(I-s\varphi(b))^{-1}
\varphi(c)(I-t\varphi(d))^{-1}\right)\,.\endalign
$$ The relation (i) now follows from Corollary 4.2 (ii) and the Neumann
series for $(I-x)^{-1}$ applied 4 times.

When computing coefficients for higher powers of $u$ one can see that
much more is true, but the combinatorics becomes very complicated.

The relation (ii) is a lot easier to prove, let
 $$x=\pmatrix
0&\;a\\
b&\;0\endpmatrix\,,
 $$
 then
 $$\align
\varphi_2(x^{2i})&=\pmatrix
\varphi((ab)^i)&\quad 0\\
0&\quad \varphi((ba)^i)\endpmatrix\\
&\\
\varphi_2(x)^{2i}&=\pmatrix
(\varphi(a)\varphi(b))^i &\quad 0\\
0&\quad (\varphi(b)\varphi(a))^i\endpmatrix\,,
\endalign
 $$
but
 $$\align
\tr_n\left((\varphi(a)\varphi(b))^i\right)&
=\tr_n\left(\varphi(a)(\varphi(b)\varphi(a))^{i-1}\varphi(b)\right)\\
&=\tr_n\left((\varphi(b)\varphi(a))^i\right)\endalign
 $$
 and
 $$\align
\tr_n\left(\varphi((ab)^i)\right)&=
\tr_n\left(\varphi(a)\varphi((ba)^{i-1}b)\right)\\
&=\tr_n\left(\varphi((ba)^i)\right)\,,\endalign
 $$
 so
 $$\align 
\tr_n\left((\varphi(a)\varphi(b))^i\right)
&=\frac12\tr_{2n}\left(\varphi_2(x)^{2i}\right)\\
&=\frac12\tr_{2n}\left(\varphi_2(x^{2i})\right)\\
&=\tr_n\left(\varphi(ab)^i)\right)\endalign
 $$
 and the proposition follows.\hfill$\square$
\enddemo

\proclaim{4.9 Corollary} Let $\varphi$ be a linear unital
 mapping of a unital matrix algebra $\Cal A$
onto a matrix algebra $\Cal B$ contained in $M_n(\Bbb C)$.

If $\varphi$ is $2$-invertibility preserving, then it is a
homomorphism modulo the Jacobson radical.
\endproclaim

\demo{Proof} By Proposition 4.8 we have for $a,b,c\in\Cal A$
 $$\align
&\tr_n\left((\varphi(a)\varphi(b)-\varphi(ab))\varphi(c)\right)\\
=&\tr_n\left(\varphi(abc)-\varphi(abc)\right)=0\,.\endalign
 $$  Since  $\varphi$ is surjective the corollary follows from Lemma
2.1.\hfill$\square$
\enddemo

\Refs{\phantom{}}

\ref\key1\by Aupetit, B.\paper Une generalisation du th\'eor\`eme de
Gleason-Kahane-\`Zelazko pour les algebres de Banach\jour Pac\. J\.
Math.\vol85\yr1979)\pages11--17
\endref

\ref\key2\by Behrens, E.-A.\book Ring theory\publ Academic
Press\yr1972
\endref

\ref\key3\by Christensen, E.\paper Two generalizations of the
Gleason-Kahane-\`Zelazko Theorem\jour Pac\. J\. Math.\toappear
\endref

\ref\key4\by Friedland, S.\paper Simultaneous similarity of
matrices\jour Adv\. in Math.\vol50\yr1983\pages189--265
\endref

\ref\key5\by Gantmacher, F\. H.\book Matrizenrechnung\publ VEB
Deutsche Verlag der Wissenschaften\publaddr Berlin\yr1958
\endref

\ref\key6\by Gleason, A\. M.\paper A characterization of maximal
ideals\jour J\. Analyse Math.\vol19\yr1967\pages171--172
\endref

\ref\key7\by Humphreys J\. E.\paper Introduction to Lie algebras and
representation theory\jour Grad\. Text in Math.\vol9 {\rm Springer}\yr1972
\endref

\ref\key8\by Kahane, J\. P\. and \`Zelazko, W.\paper A
characterization of maximal ideals in commutative Banach
algebras\jour Studia Math.\vol29\yr1968\pages339--343
\endref

\ref\key9\by Kaplansky, I.\paper Algebraic and analytic aspects of
operator algebras\jour CBMS Regional conference series in
math.\vol1\yr1970
\endref

\ref\key10\by Laffey, T\. J.\paper Simultaneous quasidiagonalization
of a pair of $3\times3$ complex matrices\jour Rev\. Roum\. Pures et
Appl.\vol23\yr1978\pages1047--1052
\endref

\ref\key11\by Laffey, T\. J.\paper Simultaneous triangularization of
matrices --- low rank cases and the nonderogatory case\jour Linear
and Multilinear Algebra\vol6\yr1978/79\pages269--305
\endref

\ref\key12\by Marcus, M\. and Ming, H.\book A survey of matrix theory
and matrix inequalities\publ Allyn and Bacon, Inc.\publaddr
Boston\yr1964
\endref

\ref\key 13\by Marcus, M\. and Purves, R.\paper Linear transformation
on algebras of matrices: {\rm The invariance of the elementary
Symmetric function}\jour Canad\. J\. Math.\vol11\yr1959\pages 383--396
\endref

\ref\key14\by McCoy, N. H. \paper On the characteristic roots of matrix polynomials \jour Bull\. Amer\. Math\.
Soc.\vol42 \yr1936\pages 592--600
\endref

\ref\key15\by Motzkin, T\. S. and Taussky, O.\paper Pairs of matrices
with property $L$\jour Trans\. Amer\. Math\.
Soc.\vol73\yr1952\pages108--114
\endref

\ref\key16\by Taussky, O.\paper Some results concerning the
transition from the $L$- to the $P$-property for pairs of finite
matrices\jour J\. Algebra\vol20\yr1972\pages271--283
\endref

\ref\key17\by Taussky, O.\paper Some results concerning the
transition from the $L$- to the $P$-property for pairs of finite
matrices {\rm II}\jour Linear and Multilinear
Algebra\vol2\yr1974\pages195--202
\endref

\ref\key18\by Robinson, N\. F. and Smiley, M\. F.\paper A
generalization of the property $L$ of Motzkin and Taussky\jour Linear
Algebra and Appl.\vol4\yr1971\pages323--328
\endref

\ref\key19\by Wielandt, H.\paper Lineare Scharen von Matrizen mit
reelle Eigenwerten\jour Math\. Zeit.\vol53\yr1950\pages219--225
\endref
\endRefs

\end
\bye